\documentclass[preprint,12pt]{elsarticle}

\usepackage{amssymb}
\usepackage{amsthm}

\usepackage{multicol}
\usepackage{mathrsfs}
\usepackage{amsfonts}
\usepackage{CJK,fancybox,color,graphicx,amsmath,amsthm,amssymb,enumerate}

\newtheorem{thm}{Theorem}[section]

\newtheorem{lem}[]{Lemma}[section]

\theoremstyle{definition}
\newtheorem{defn}[thm]{Definition}
\theoremstyle{remark}
\newtheorem{rem}[]{Remark}[section]
\numberwithin{equation}{section}

\newcommand{\D}{\displaystyle}
\newcommand{\DF}[2]{\frac{\D#1}{\D#2}}

\begin{document}

\begin{frontmatter}



\title{\bf The Riemann solution to the Chaplygin pressure Aw-Rascle
model with Coulomb-like friction and its vanishing pressure limit}

\author{Qingling Zhang}

\address{School of Mathematics and Computer Sciences, Jianghan University, Wuhan 430056, PR China }
\ead{zhangqingling2002@163.com}
\begin{abstract}

 The Riemann solution to the Chaplygin pressure Aw-Rascle model with
Coulomb-like friction is constructed explicitly and its vanishing
pressure limit is analyzed precisely. It is shown that the delta
shock wave appears in the Riemann solutions in some certain
situations. The generalized Rankine-Hugoniot conditions of the delta
shock wave are established and the exact position, propagation speed
and strength of the delta shock wave are given explicitly, which
enables us to see the influence of the Coulomb-like friction on the
Riemann solution to the Chaplygin pressure Aw-Rascle model clearly.
It is shown that the Coulomb-like friction term makes contact
discontinuities and delta shock waves bend into parabolic shapes and
the Riemann solutions are not self-similar anymore. Finally, the
occurrence mechanism on the phenomenon of concentration and
cavitation and the formation of delta shock wave and vacuum in the
process of vanishing pressure limit are analyzed and identified in
detail. Moreover, we show the Riemann solutions to the
nonhomogeneous Chaplygin pressure Aw-Rascle model converge to the
Riemann solutions to the transportation equations with the same
source term as the pressure vanishes. These two results generalize
those obtained in \cite{Chen-Liu1,Sheng-Wang-Yin} for homogeneous
equations to nonhomogeneous equations and are also applicable to the nonsymmetric system of
Keyfitz-Kranzer type with the same Chaplygin pressure and Coulomb-like friction.

\end{abstract}


\end{frontmatter}



\section{Introduction}
\setcounter{equation}{0}

\ \  In this paper, we are mainly concerned with the Riemann problem
for the Chaplygin pressure Aw-Rascle model with Coulomb-like
friction
\begin{equation}\label{1.1}
\left\{\begin{array}{ll}
\rho_t+(\rho u)_x=0,\\
(\rho( u+P))_t+(\rho u( u+P))_x=\beta\rho,
 \end{array}\right.
\end{equation}
 with Riemann initial data
\begin{equation}\label{1.2}
(\rho,u)(x,0)=\left\{\begin{array}{ll}
(\rho_-,u_-),\ \ x<0,\\
(\rho_ +,u_+),\ \ x>0.
\end{array}\right.
\end{equation}
where $\rho_{\pm}$ and $u_{\pm}$ are all given constants. In
(\ref{1.1}), the state variable $\rho>0$ and $u\geq0$ denote the
traffic density and velocity, respectively, $\beta$ is a frictional
constant, and the pressure $P$ is given by the state equation
\begin{equation}\label{1.3}
P=-\frac{A}{\rho}, \ \ A>0,
\end{equation}
which was introduced by Chaplygin \cite{Chaplygin} and
Tsien\cite{Tsien} as a suitable mathematical approximation for
calculating the lifting force on a wing of an airplane in
aerodynamics.

The Euler system with state equation (\ref{1.3}) is the classical
Chaplygin gas equations which has been advertised as a possible
model for dark energy of the universe
\cite{Bilic-Tupper-Viollier,Gorini-Kamenshchik-Moschella-Pasquier}
and have been extensively investigated recently
\cite{Brenier,Guo-Sheng-Zhang,Serre,Wang-Zhang} etc.  The
generalized Chaplygin gas model has also attracted intensive
attention such as in
\cite{Bilic-Tupper-Viollier,Setare,Sheng-Wang-Yin,Wang}. It can be
used to describe the dark matter and dark energy in the unified form
through exotic background fluid whose state  equation is given by
$P=-\frac{A}{\rho^{\alpha}},0<\alpha<1,A>0$. The modified Chalygin
gas was proposed by Benaoum in 2002 \cite{Benaoum} to describe the
current accelerated expansion of the universe, whose equation of
state is given by $P=A\rho-\frac{B}{\rho^{\alpha}},\ \
0<\alpha\leq1, A,B>0$. Compared with the Chalygin gas or the
generalized Chalygin gas, the model for the modified Chalygin gas
can describe the universe to a large extent.

If $\beta=0$, then the system (\ref{1.1}) becomes the Chaplygin
pressure Aw-Rascle model which was recently introduced by Pan and
Han \cite{Pan-Han}, in which delta-shocks appear in the Riemann
solutions, which may be used to explain the serious traffic jam.
Sheng and Zeng \cite{Sheng-Zeng} considered its Riemann problem with
delta initial data. With these results, similar to
\cite{Huang-Wang,Wang-Ding,Wang-Huang-Ding}, we recently solved the
Cauchy problem of it in \cite{Zhang2} by generalized potential
method. If $\beta=0$ and $P=\rho^{\gamma},\gamma>0$, then the system
(\ref{1.1}) becomes the classical Aw-Rascle model of traffic flow
proposed by Aw and Rascle \cite{Aw-Rascle} in 2000 to remedy the
deficiencies of second order models of car traffic pointed out by
Daganzo \cite{Daganzo} and had also been independently derived by
Zhang \cite{Zhang3}. Since then, it had received extensive attention
\cite{Greenberg,Lebacque-Mannnar-Salem,Shen-Sun,Sun}. Recently, the
Riemann problem for the Aw-Rascle model with generalized Chaplygin
pressure was also considered by Guo in \cite{Guo} in which the delta-shock
also appears. Cheng and Yang \cite{Cheng-Yang} considered the
Riemann problem for the Aw-Rascle model with modified Chaplygin
pressure $P=A\rho-\frac{B}{\rho},A,B>0$ and analyzed the limit of
its Riemann solutions with the pressure approaching Chaplygin
pressure.

In fact, if $\beta=0, P=\rho^{\gamma},\gamma>0 $ and let $u=w-P$ , then the classical Aw-Rascle model can be written as the nonsymmetric
system of Keyfitz-Kranzer type as follows:
\begin{equation}\label{1.4}
\left\{\begin{array}{ll}
\rho_t+(\rho (w-P))_x=0,\\
(\rho w)_t+(\rho w( w-P))_x=0,
 \end{array}\right.
\end{equation}
Recently, Lu \cite{Lu} studied the existence of global entropy solution to general system of
Keyfitz-kraner type (\ref{1.4}) with state equation $P=P(\rho)$ satisfying some conditions.
In 2013, Cheng \cite{Cheng1,Cheng2} considered the Riemann problem of (\ref{1.4}) with different choice of state
equation of $P$, such as $P$ taken as the Chaplygin pressure, the generalized Chaplygin pressure and the
modified Chaplygin pressure, etc, which showed that the Riemann solutions to (\ref{1.4}) with Chaplygin pressure and
generalized Chaplygin pressure were very similar to that of the Aw-Rascle model with the
corresponging pressure.

If $\beta=0$ and $P=0$, then the system (\ref{1.1}) becomes the
so-called zero pressure flow (transportation equations). It is well
known that the delta-shock wave also appears in the Riemann
solutions to the zero pressure flow which has been widely studied
such as
\cite{Bouchut,E-Rykov-Sinai,Huang-Wang,Li-Zhang-Yang,Sheng-Zhang,Wang-Ding,Wang-Huang-Ding}.
Recently, Shen \cite{Shen1} considered (\ref{1.1}) with $P=0$ and
solved the Riemann problem and the generalized Riemann problem for
the transportation equations Coulomb-like friction. Delta-shock is a
very interesting topic in the theory of conservation laws. It is a
generalizations of an ordinary shock. Speaking informally, it is a
kind of discontinuity, on which at least one of the variables may be
develop an extreme concentration in the form of a weighted Dirac
delta function with the discontinuity as its support. From the
physical point of view, it represents the process of the
concentration of the mass. For related research of delta-shock
waves, we refer readers to papers
\cite{Korchinski,Li-Zhang-Yang,Sheng-Zhang,Tan-Zhang,Tan-Zhang-Zheng}
and the references cited therein for more details.

 From the above discussions, one can see that the Riemann problem for the
Aw-Rascle model with various kinds of pressure but without source
term (namely $\beta=0$) has been well investgated. Hence, it is
natural to expect the study of it with a source term, such as
damping, friction and relaxation effect. In the present paper, we
want to deal with the Riemann problem for the Chaplygin pressure
Aw-Rascle model with Coulomb-like friction which was proposed by
Savage and Hutter in 1989 \cite{Savage-Hutter} to describe granular
flow behavior. For research on other models with Coulomb-like
friction, one can see \cite{Shen1,Shen2,Sun1}.

In this paper, we are interested in how the delta-shock solution of
the Chaplygin pressure Aw-Rascle model with Coulomb-like friction
develops under the influence of the Coulomb-like friction. The
advantage of this kind source term is in that (\ref{1.1}) can be
written in a conservative form such that exact solutions to the
Riemann problem (\ref{1.1}) and (\ref{1.2}) can be constructed
explicitly. We shall see that the Riemann solutions to  (\ref{1.1})
and (\ref{1.2}) are not self-similar any more, in which the state
variable $u$ varies linearly along with the time $t$ under the
influence of the Coulomb-like friction. In other words, the state
variable $u-\beta t$ remains unchanged in the left, intermediate and
right states. In some situations, the delta-shock wave appears in
the Riemann solutions to  (\ref{1.1}) and (\ref{1.2}). In order to
describe the delta-shock wave, the generalized Rankine-Hugoniot
conditions are derived and the exact position, propagation speed and
strength of the delta shock wave are obtained completely. It is
shown that the Coulomb-like friction term make contact
discontinuities and delta shock waves bend into parabolic shapes for
the Riemann solutions.

 Finally, the occurrence mechanism on the
phenomenon of concentration and cavitation and the formation of
delta shock wave and vacuum in the process of vanishing pressure
limit of Riemann solutions to the nonhomogeneous Chaplygin pressure
Aw-Rascle model are analyzed and identified in detail, from
which we find that there is something different from polytropic gas
in \cite{Chen-Liu1} but similar to generalized Chaplygin gas in
\cite{Sheng-Wang-Yin} about the formation of the delta shock wave.
Moreover, we show the Riemann solutions to the nonhomogeneous
Chaplygin pressure Aw-Rascle model converge to the Riemann solutions
for the transportation equations with the same source term as the
pressure vanishes. These two results generalize those obtained in
\cite{Chen-Liu1,Sheng-Wang-Yin} for homogeneous equations to
nonhomogeneous equations. Since the configuration of the Riemann solution to (\ref{1.4}) with Chaplygin pressure
is very similar to that of (\ref{1.1}) with $\beta=0$ (see \cite{Cheng1}), we can obtain similar
results for the nonsymmetric system of Keyfitz-Kranzer type (\ref{1.4}) with the same Chaplygin pressure and Coulomb-like friction.

 In fact, it was shown in \cite{Guo} that the delta-shock also appears in
the Riemann solutions to the Aw-Rascle model with generalized
Chaplygin pressure.  It should be remarkable that a significant
mathematical difference between the Aw-Rascle model with generalized
Chaplygin pressure and with Chaplygin pressure for the reason that
there is one characteri tic field genuinely nonlinear for the
former, whose elementary waves admit not only contact
discontinuities, but also rarefaction waves and shock waves, while
the two characteristic fields are all linearly degenerate for the
latter, whose elementary waves admit only contact discontinuities.
To investigate how the Coulomb-like friction affects the rarefaction
 aves, shock waves and the delta shock waves and the occurrence mechanism of the delta
shock waves in the process of pressure decreasing, we will study the
Riemann problem for the generalized Chaplygin pressure Aw-Rascle
model with Coulomb-like friction and its vanishing pressure limit, whose results will
also be applicable to the nonsymmetric system of Keyfitz-Kranzer type (\ref{1.4}) with the same
pressure and Coulomb-like friction.

This paper is organized as follows. In Section 2, the system
(\ref{1.1}) is reformulated into a conservative form and then some
general properties of the conservative form are obtained. Then, the
exact solution to the Riemann problem  for the  conservative form
are constructed explicitly, which involves the delta shock wave.
Furthermore, the generalized Rankine-Hugoniot conditions are
established and the exact position, propagation speed and strength
of the delta shock wave are given explicitly. In Section 3, the
generalized Rankine-Hugoniot conditions and the exact Riemann
solutions to (\ref{1.1}) and (\ref{1.2}) are also given.
Furthermore, it is proven rigorously that the delta-shock wave is
indeed a week solution to the Riemann problem (\ref{1.1}) and
(\ref{1.2}) in the sense of distributions. In Section 4, we analyze
the formation of delta shock waves and vacuum states in the Riemann
solutions to (\ref{1.1}) and (\ref{1.2}) in the vanishing pressure
limit and show that the Riemann solutions converge to the
corresponding ones of the transportation equaitons with the same
source term as the pressure vanishes. Finally, conclusions and
discussions are carried out in Section 5.

\section{Riemann problem for a modified conservative system}
In this section, we are devoted to the study of the Riemann problem
for a conservative system (\ref{1.1}) in detail. Let us introduce
the new velocity $v(x,t)=u(x,t)-\beta t$, then the system
(\ref{1.1}) can be reformulated in a conservative form as follows:

\begin{equation}\label{2.1}
\left\{\begin{array}{ll}
\rho_t+(\rho (v+\beta t))_x=0,\\
(\rho( v+P))_t+(\rho (v +P)(v+\beta t))_x=0.
 \end{array}\right.
\end{equation}
In fact, the change of variable was introduced by Faccanoni and
Mangeney \cite{Faccanoni-Mangeney} to study the shock and
rarefaction waves of the Riemann problem for the shallow water
equations with a with Coulomb-like friction term. Here, we use this
transformation to study the delta shock wave for the system
(\ref{1.1}) which is a fully linearly degenerate system.

Now we want to deal with the Riemann problem for the conservative
system (\ref{2.1}) with the same
 Riemann initial data (\ref{1.2}) as follows:
\begin{equation}\label{2.2}
(\rho,v)(x,0)=\left\{\begin{array}{ll}
(\rho_-,u_-),\ \ x<0,\\
(\rho_ +,u_+),\ \ x>0.
\end{array}\right.
\end{equation}
We shall see hereafter that the Riemann solutions to (\ref{1.1}) and
(\ref{1.2}) can be obtained immediately from the Riemann solutions
to (\ref{2.1}) and (\ref{2.2}) by using the transformation of state
variables $(\rho,u)(x,t)=(v+\beta t)(x,t)$.

The system (\ref{2.1}) can be rewritten in the quasi-linear form
\begin{equation}\label{2.3}
\left(\begin{array}{lll}
1 & 0 \\
v & \rho
\end{array}\right)
\left(\begin{array}{lll}
\rho  \\
v
\end{array}\right)_{t}
+\left(\begin{array}{lll}
v+\beta t & \rho  \\
v(v+\beta t) &\rho(2v+\beta t+P)
\end{array}\right)
\left(\begin{array}{lll}
\rho  \\
v
\end{array}\right)_{x}
=\left(\begin{array}{lll}
0  \\
0
\end{array}\right).
\end{equation}
It can be derived directly from (\ref{2.3}) that the conservative
system (\ref{2.1}) has two eigenvalues
$$\lambda_1(\rho,v)=v+\beta t-\frac{A}{\rho},\ \ \lambda_2(\rho,v)=v+\beta t,$$
whose corresponding right eigenvectors can be expressed respectively
by
$$r_1=(\rho,-\frac{A}{\rho})^T,\ \ r_2=(1,0)^T.$$
So (\ref{2.1}) is strictly hyperbolic for $\rho>0$. Moreover,
$\bigtriangledown\lambda_i\cdot r_i=0$, $i=1,2$, which implies that
$\lambda_1$ and $\lambda_2$ are both linearly degenerate and the
associated waves are both contact discontinuities denoted by $J$,
see \cite{Smoller}.

 We should take notice the fact that the parameter $t$ only appears
 in the flux functions in the conservative
system  (\ref{2.1}), such that the Ranking-Hugoniot conditions can
be derived in a standard method as in \cite{Smoller}. For a bounded
discontinuity at $x=x(t)$, let us denote $\sigma(t)=x'(t)$, then the
Rankine-Hugoniot conditions for the conservative system
(\ref{2.1})can be expressed as
\begin{equation}\label{2.4}
\left\{\begin{array}{ll}
-\sigma(t)\rho+[\rho (v+\beta t)]=0,\\
-\sigma(t)[\rho (v+P)]+[\rho(v+P) (v+\beta t)]=0,
\end{array}\right.
\end{equation}
where $[\rho]=\rho_{r}-\rho_l$ with
$\rho_{l}=\rho(x(t)-0,t)$,\ $\rho_{r}=\rho(x(t)+0,t)$, in which
$[\rho]$ denote the jump of $\rho$ across the discontinuity, etc. It
is clear that the propagation speed of the discontinuity depends on
the parameter $t$, which is obviously different from classical
hyperbolic conservation laws.

If $\sigma(t)\neq0$, then it follows from (\ref{2.4}) that

\begin{equation}\label{2.5}
\rho_{r}\rho_{l}(v_{r}-v_{l})((v_{r}-\frac{A}{\rho_{r}})-(v_{l}-\frac{A}{\rho_{l}}))=0,
\end{equation}
from which we have $v_{r}=v_{l}$ or
$v_{r}-\frac{A}{\rho_{r}}=v_{l}-\frac{A}{\rho_{l}}$.

 Thus, the two states $(\rho_{r},v_{r})$ and $(\rho_{l},v_{l})$ can
be connected by a 1-contact discontinuity if and only if
\begin{equation}\label{2.6}
J_{1}:\ \  \ \sigma(t)=v_{r}+\beta t-\frac{A}{\rho_{r}}=v_{l}+\beta
t-\frac{A}{\rho_{l}},
\end{equation}
and can be connected by a 2-contact discontinuity if and only if
\begin{equation}\label{2.7}
J_{2}:\ \  \ \sigma(t)=v_{r}+\beta t=v_{l}+\beta t.
\end{equation}


\unitlength 0.9mm 
\linethickness{0.4pt}
\ifx\plotpoint\undefined\newsavebox{\plotpoint}\fi 
\begin{picture}(200,69)(-15,0)
\put(95.75,13.75){\vector(1,0){.07}}
\put(15.75,13.75){\line(1,0){80}}
\put(15.75,69){\vector(0,1){.07}}
\put(15.75,14){\line(0,1){55}}
\bezier{50}(35.25,67.75)(35.25,23.25)(35.25,13.75)
\put(65,68){\line(0,-1){53.75}}
\qbezier(37.5,65.5)(46,23.25)(88.5,21)
\put(12.5,67.5){$\rho$}
\put(99,13.75){$v$}
\put(15.5,9.75){0}
\put(26,41.25){I\!I\!I}
\put(55,46.25){I\!I}
\put(76.5,45.5){I}
\put(44.75,53.75){$J_1$}
\put(69,55.25){$J_2$}
\put(31,50.75){$S$}
\put(30,10){$u_--\frac{A}{\rho_-}$}
\put(63.5,11.5){$u_-$}
\put(69,28.5){($\rho_-$,$u_-$)}
\put(60.75,5.75){\makebox(0,0)[cc]
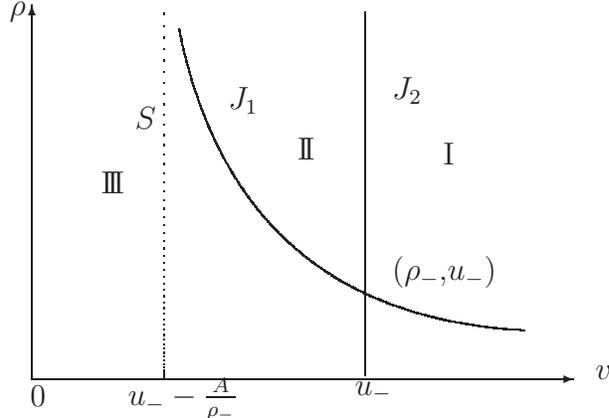
{Fig.1 the $(\rho,v)$ phase plane for the conservative system (2.1).}}
\end{picture}

In the $(\rho,v)$ phase plane, for the given state $(\rho_-,u_-)$,
it follow from (\ref{2.6}) that the sets of states connected on the
right consists of the 1-contact discontinuity curve
$J_{1}(\rho_-,u_-)$ satisfying
$v-\frac{A}{\rho}=u_{-}-\frac{A}{\rho_{-}}$, which has two
asymptotes $S: v=u_{-}-\frac{A}{\rho_{-}}$ and $\rho=0$. Similarly,
it follow from (\ref{2.7}) that the sets of states connected on the
right consists of the 2-contact discontinuity curve
$J_{2}(\rho_-,u_-)$ satisfying $v=u_{-}$.  In the $(\rho,v)$ phase
plane with $\rho>0,v\geq0$, let us draw Fig.1 to depict these curves
together which divide the $(\rho,v)$ phase plane into three parts I,
I\!I and I\!I\!I,where
$${\rm I}=\{(\rho,v)|v\geq u_-\},$$
$${\rm I\!I}=\{(\rho,v)|u_{-}-\frac{A}{\rho_{-}}<v< u_-\},$$
$${\rm I\!I\!I}=\{(\rho,v)|v\leq u_--\frac{A}{\rho_{-}}\}.$$

 When
$(\rho_+,u_+)\in$ I $\cup$ I\!I, namely
$u_{+}>u_{-}-\frac{A}{\rho_{-}}$, the Riemann solutions consists of
 two contact discontinuity $J_{1}$ and $J_{2}$ with the intermediate constant state
$(\rho_\ast,v_\ast)$ between them besides constant states
$(\rho_-,u_-)$ and $(\rho_+,u_+)$, where
\begin{equation}\label{2.8}
\left\{\begin{array}{ll}
v_{\ast}-\frac{A}{\rho_{\ast}}=u_--\frac{A}{\rho_-},\\
u_+=v_{\ast}.
\end{array}\right.
\end{equation}
which immediately leas to

\begin{equation}\label{2.9}
(\frac{A}{\rho_\ast},v_\ast)=(u_+-u_-+\frac{A}{\rho_-},u_+).
\end{equation}
The propagation speed of $J_{1}$ and $J_{2}$ are given by
$\sigma_{1}(t)=u_--\frac{A}{\rho_-}+\beta t$ and
$\sigma_{2}(t)=u_++\beta t$, respectively.

On the other hand, when $0\leq(\rho_+,u_+)\in$ I\!I\!I, namely
$u_{+}\leq u_{-}-\frac{A}{\rho_{-}}$, then the characteristic curves for
the Riemann problem (\ref{2.1}) and (\ref{2.2}) overlap in a domain
$\Omega$ such that singularity will happen in $\Omega$. For
completeness, we simply compute the characteristic curves emitting
from the origin $(0,0)$ which are determined by
$$\frac{dx^{\pm}_{i}(t)}{dt}=\lambda_{i}(\rho_{\pm},u_{\pm}).$$
Thus, we have
$$x^{-}_{1}(t)=(u_--\frac{A}{\rho_-})t+\frac{1}{2}\beta t^{2},\ \ x^{+}_{1}(t)=(u_+-\frac{A}{\rho_+})t+\frac{1}{2}\beta t^{2},$$
$$x^{-}_{2}(t)=u_-t+\frac{1}{2}\beta t^{2},\ \ x^{+}_{2}(t)=u_+t+\frac{1}{2}\beta t^{2}.$$
Let us draw Fig.2 to explain this phenomenon in detail. In fact, the
Cauchy problem for the Chaplygin pressure Aw-Rascle model has been
well investigated by us \cite{Zhang2} recengly by using the
generalized characteristic method.

\unitlength 1mm 
\linethickness{0.4pt}
\ifx\plotpoint\undefined\newsavebox{\plotpoint}\fi 
\begin{picture}(106,66)(-10,0)
\put(102,10.25){\vector(1,0){.07}}
\put(24.75,10.25){\line(1,0){77.25}}
\put(19,65.75){\vector(0,1){.07}}
\put(19,12.25){\line(0,1){53.5}}
\qbezier(56.25,10)(77.88,33.25)(100.25,40.5)
\qbezier(56.25,10)(81.25,28.38)(102.25,31)
\qbezier(56.25,10)(35.25,28.63)(47.75,61.75)
\put(15,63.25){$t$}
\put(106,10){$x$}
\put(48.5,62.5){$x_1^+(t)$}
\put(72,62.5){$x_2^+(t)$}
\put(83.25,54.5){$\delta \!S$}
\put(104,42.75){$x_1^-(t)$}
\put(104.5,31.25){$x_2^-(t)$}
\put(58.5,28){$<$}
\put(67.5,22){$\vee$}
\put(64,30){$\Omega$}
\put(56.25,6.5){0}
\put(70,5){\makebox(0,0)[cc]
{Fig.2 The characteristic analysis of delta shock wave for the }}
\put(65,0.5){\makebox(0,0)[cc]
{Riemann problem (2.1) and (2.2) when $u_{+} < u_{-}-\frac{A}{\rho_{-}}$.}}

\qbezier(56.25,10)(53.5,30.25)(73.5,59.5)
\bezier{50}(56.25,10)(65.63,35.88)(86.25,51.5)
\qbezier(58.25,29.75)(67,30)(69,23)
\end{picture}

The formation of singularity for the solution to Riemann problem
(\ref{2.1}) and (\ref{2.2}) is due to the overlap of linearly
degenerate characteristics. Thus, the nonclassical situation appears
for some certain initial data where the Cauchy problem usually does
not own a weak $L^{\infty}$-solution. In order to solve the Riemann
problem (\ref{2.1}) and (\ref{2.2}) in the framework of nonclassical
solutions, motivated by \cite{Pan-Han}, a solution containing a
weighted $\delta$-measure supported on a curve should be introduced.

\begin{defn}\label{defn:2.1}

To define the measure solutions, a two-dimensional weighted
$\delta$-measure $p(s)\delta_S$ supported on a smooth curve
$S=\{(x(s),t(s)):a<s<b\}$ can be defined as
\begin{equation}\label{2.10}
\langle
p(s)\delta_S,\psi(x(s),t(s))\rangle=\int_a^bp(s)\psi(x(s),t(s))\sqrt{{x'(s)}^2+{t'(s)}^2}ds,
\end{equation}
for any $\psi\in C_0^\infty(R\times R_{+})$.
\end{defn}

For convenience, we usually select the parameter $s=t$ and use
$w(t)=\sqrt{1+{x'(t)}^2}p(t)$ to denote the strength of delta shock
wave from now on. In what follows, let us provide the definition of
delta shock wave solution to the Riemann problem (\ref{2.1}) and
(\ref{2.2}) in the framework introduced by Danilov and Shelkovich
\cite{Danilvo-Shelkovich1,Danilvo-Shelkovich2} and developed by
Kalisch and Mitrovic \cite{Kalisch-Mitrovic1,Kalisch-Mitrovic2}.

Let us suppose that $\Gamma=\{\gamma_{i}\mid i\in I\}$ is a graph in
the upper half plane $\{(x,t)\mid x\in R, t\in [0,+\infty)\}$, which
contains Lipschtitz continuous arcs $\gamma_{i}$ where $ i\in I$ and
$I$ is a finite index set. Let $I_{0}$ be s subset of $I$ which
contains all indices of arcs starting from the $x$-axis. Let us use
$\Gamma_{0}=\{x_{j}^{0}\mid j\in I_{0}\}$ to denote the set of
initial points of the arcs $\gamma_{j}$ with $j\in I_{0}$. Then ,one
can define the solutions in the sense of distributions to Cauchy
problem for the conservative system (\ref{2.1}) with delta measure
initial data below.

\begin{defn}\label{defn:2.2}

Let $(\rho,v)$ be a pair of distributions where $\rho$ is
represented in the form
\begin{equation}\label{2.11}
\rho(x,t)=\hat{\rho}(x,t)+w(x,t)\delta(\Gamma),
\end{equation}
in which $\hat{\rho},v\in L^{\infty}(R\times R_{+})$ and the
singular part is defined by
\begin{equation}\label{2.12}
w(x,t)\delta(\Gamma)=\sum\limits_{i\in
I}w_{i}(x,t)\delta(\gamma_{i}).
\end{equation}
Let us consider the delta shock wave type initial data of the form
\begin{equation}\label{2.13}
(\rho,v)(x,0)=(\hat{\rho}(x)+\sum\limits_{j\in
I_{0}}w_{j}(x_{j}^{0},0)\delta(x-x_{j}^{0}),v_{0}(x)).
\end{equation}
in which $\hat{\rho}_{0}(x),v_{0}(x)\in L^{\infty}(R)$, then the
above pair of distributions $(\rho,v)$ are called as a generalized
delta shock wave solution to the conservative system (\ref{2.1})
with initial data (\ref{2.19}) if the following integral identities
\begin{eqnarray}
&&\int_{R_{+}}\int_{R}(\hat{\rho}\psi_{t}+\hat{\rho}(v+\beta
t)\psi_{x}+)dxdt+\sum\limits_{i\in
I}\int_{\gamma_{i}}w_{i}(x,t)\frac{\partial\psi(x,t)}{\partial l}dl\nonumber\\
&&+\int_{R}\hat{\rho}_{0}(x)\psi(x,0)dx+\sum\limits_{k\in
I_{0}}w_{k}(x_{k}^{0},0)\psi(x_{k}^{0},0)=0\label{2.14},
\end{eqnarray}
\begin{eqnarray}
&&\int_{R_{+}}\int_{R}(\hat{\rho}(v+P)\psi_{t}+\hat{\rho}(v+P)(v+\beta
t)\psi_{x}dxdt+)+\sum\limits_{i\in
I}\int_{\gamma_{i}}w_{i}(x,t)v_{\delta}(x,t)\frac{\partial\psi(x,t)}{\partial l}dl\nonumber\\
&&+\int_{R}\hat{\rho}_{0}(x)v_{0}(x)\psi(x,0)dx+\sum\limits_{k\in
I_{0}}w_{k}(x_{k}^{0},0)v_{\delta}(x_{k}^{0},0)(\psi(x_{k}^{0},0)=0\label{2.15},
\end{eqnarray}
hold for any test function  $\psi\in C_c^\infty(R\times R_{+})$, in
which $\frac{\partial\psi(x,t)}{\partial l}$ stands for the
tangential derivative of a function $\psi$ on the graph $\gamma_{i}$
and $\int_{\gamma_{i}}$ is the line integral along the arc
$\gamma_{i}$.
\end{defn}

The above-defined singular solution should be understood in the
sense of weak asymptotic solutions. More precisely, let
$f_{\epsilon}(x)\in D'(R)$ be a family of distributions depending on
$\epsilon\in(0,1)$, then we have $f_{\epsilon}(x)=o_{D'}(1)$ if the
estimate $\langle f_{\epsilon,\psi}\rangle=o(1)$ as
$\epsilon\rightarrow0$ holds for any $\psi\in D(R)$. Then, the
family of pairs of functions $(\rho_{\epsilon},v_{\epsilon})$ is
called a weak asymptotic solution of the Cauchy problem (\ref{2.1})
and (\ref{2.12}) if the limit of $\epsilon\rightarrow0$ of
$(\rho_{\epsilon},v_{\epsilon})$ is a pair of distributions for
every fixed $t\in R_{+}$, where

\begin{equation}\label{2.16}
\left\{\begin{array}{ll}
(\rho_{\epsilon})_t+(\rho_{\epsilon} (v_{\epsilon}+\beta t))_x=o_{D'}(1),\\
(\rho_{\epsilon}( v_{\epsilon}+P_{\epsilon}))_t+(\rho_{\epsilon}
(v_{\epsilon} +P_{\epsilon})(v_{\epsilon}+\beta t))_x=o_{D'}(1),
 \end{array}\right.
\end{equation}
and
\begin{equation}\label{2.17}
\rho_{\epsilon}|_{t=0}-\rho(x,0)=o_{D'}(1),\ \
v_{\epsilon}|_{t=0}-v(x,0)=o_{D'}(1)\ \  as \epsilon\rightarrow0.
\end{equation}

It can be seen from \cite{Danilvo-Shelkovich1,Danilvo-Shelkovich2}
that the limit $(\rho_{\epsilon},v_{\epsilon})$ as
$\epsilon\rightarrow0$ can be understood in Definition 2.2. The weak
asymptotic solution is constructed such that the terms that do not
have a distributional limit cancel in the limit
$\epsilon\rightarrow0$ and the problem about multiplication of
distributions is automatically eliminated.

With the above definition, if $(\rho_+,u_+)\in$ I\!I\!I and $u_{+}<u_{-}-\frac{A}{\rho_{-}}$, we consider solutions of the form
\begin{equation}\label{2.18}
(\rho,v)(x,t)=\left\{\begin{array}{ll}
(\rho_-,u_-),\ \ \ \ \ \ \ \ \ \ \ \ \ \ \ \ \ \ \ \ x<x(t),\\
(w(t)\delta(x-x(t)),v_{\delta}),\ \ \  \ x=x(t),\\
(\rho_+,u_+),\ \ \ \ \ \ \ \ \ \ \ \ \ \ \ \ \ \ \ \ x>x(t),
\end{array}\right.
\end{equation}
where $x(t)$, $w(t)$ and $\sigma(t)=x'(t)$ denote respectively the
location, weight and propagation speed of the delta shock,
$v_{\delta}$ indicates the assignment of $v$ on this  delta shock
wave, and $\frac{1}{\rho}$ is equal to zero on this  delta shock
wave. In fact, the  delta  shock wave solution (\ref{2.18}) to the
the Riemann problem (\ref{2.1}) and (\ref{2.2}) is the simplest
example that the graph $\Gamma$ contains only one arc. When $u_{+}=u_{-}-\frac{A}{\rho_{-}}$,
it can be discussed similarly and we omit it.

Let us check briefly that the delta shock wave solution of the form
(\ref{2.18}) to the the Riemann problem (\ref{2.1}) and (\ref{2.2})
 satisfy the following generalized Rankine-Hugoniot conditions
\begin{equation}\label{2.19}
\left\{\begin{array}{ll}
\DF{dx(t)}{dt}=\sigma(t)=v_{\delta}+\beta t,\\[4pt]
\DF{dw(t)}{dt}=\sigma(t)[\rho]-[\rho (v+\beta t)],\\[4pt]
\DF{d(w(t)v_{\delta})}{dt}=\sigma(t)[\rho( v-\frac{A}{\rho})]-[\rho(
v-\frac{A}{\rho}) (v+\beta t)].
\end{array}\right.
\end{equation}

Let us assume that the delta shock wave curve
$\Gamma:{(x,t)|x=x(t)}$ is a smooth curve in the $(x,t)$ plane
across which $(\rho,v)$ is a jump discontinuity.Let $P$ be any point
on $\Gamma$ and let $\Omega$ be a small ball centered at the point
$P$. Then, we make a step further to assume that the intersection
point of $\Omega$ and $\Gamma$ are $P_{1}=(x(t_{1}),t_{1})$ and
$P_{2}=(x(t_{2}),t_{2})$ where $t_{1}<t_{2}$, and $\Omega_{-}$ and
$\Omega_{+}$ are the left-hand and right-hand parts of $\Omega$ cut
by $\Gamma$ respectively. Then, for any test function $\psi\in
C_c^\infty(\Omega)$, by applying the divergence theorem, we have

\begin{eqnarray}
I&=&\int\int_{\Omega}\Big(\rho( v-\frac{A}{\rho})\psi_{t}+\rho(
v-\frac{A}{\rho}) (v+\beta t)\psi_{x}\Big)dxdt\nonumber\\
&=&\int\int_{\Omega_{-}}\Big(\rho_{-}(
u_{-}-\frac{A}{\rho_{-}})\psi_{t}+\rho_{-}(u
_{-}-\frac{A}{\rho_{-}}) (u_{-}+\beta t)\psi_{x}\Big)dxdt \nonumber\\
&&\int\int_{\Omega_{+}}\Big(\rho_{+}(
u_{+}-\frac{A}{\rho_{+}})\psi_{t}+\rho_{+}(u
_{+}-\frac{A}{\rho_{+}}) (u_{+}+\beta t)\psi_{x}\Big)dxdt \nonumber\\
&&+\int_{t_{1}}^{t_{2}}w(t)\Big(v_{\delta}\psi_{t}(x(t),t)+v_{\delta}(v_{\delta}+\beta
t)\psi_{x}(x(t),t)\Big)dt\nonumber\\
&=&\int_{\partial\Omega_{-}}-\rho_{-}( u_{-}-\frac{A}{\rho_{-}})\psi
dx+\rho_{-}(u
_{-}-\frac{A}{\rho_{-}}) (u_{-}+\beta t)\psi dt \nonumber\\
&&\int_{\partial\Omega_{+}}-\rho_{+}( u_{+}-\frac{A}{\rho_{+}})\psi
dx+\rho_{+}(u_{+}-\frac{A}{\rho_{+}}) (u_{+}+\beta t)\psi dt \nonumber\\
&&+\int_{t_{1}}^{t_{2}}w(t)\Big(v_{\delta}\psi_{t}(x(t),t)+v_{\delta}(v_{\delta}+\beta
t)\psi_{x}(x(t),t)\Big)dt\nonumber
\end{eqnarray}
\begin{eqnarray}
&=&\int_{t_{1}}^{t_{2}}\Big((\rho_{+}(u_{+}-\frac{A}{\rho_{+}})-\rho_{-}(
u_{-}-\frac{A}{\rho_{-}}))\frac{dx}{dt}\nonumber\\
&&+\rho_{-}(u_{-}-\frac{A}{\rho_{-}}) (u_{-}+\beta
t)-\rho_{+}(u_{+}-\frac{A}{\rho_{+}}) (u_{+}+\beta
t)\Big)\psi(x(t),t)dt\nonumber\\
&&\int_{t_{1}}^{t_{2}}w(t)v_{\delta}d\psi(x(t),t)\nonumber.
\end{eqnarray}

Thus, one can see that the third equality in (\ref{2.19}) is
satisfied when $I$ vanishes for any $\psi\in C_c^\infty(\Omega)$. In
the same way as above, we can check that the second identity
holds.Thus, the proof is complete.

In order to ensure uniqueness, it should also satisfy an
over-compressive entropy condition for the delta shock wave as
follows:
\begin{equation}\label{2.20}
\lambda_{1}(\rho_{+},u_{+})<\lambda_{2}(\rho_{+},u_{+})<\sigma(t)<\lambda_{1}(\rho_{-},u_{-})<\lambda_{2}(\rho_{-},u_{-})£¬
\end{equation}
which enables us to have
\begin{equation}\label{2.21}
0\leq u_{+}<v_{\delta}<u_{-}-\frac{A}{\rho_{-}}.
\end{equation}

The  generalized Rankine-Hugoniot conditions (\ref{2.19}) reflect
the relationship among the location, weight and propagation speed of
delta shock wave. The entropy condition (\ref{2.20}) for the delta
shock wave is an over-compressive condition which implies that all
the characteristics on both sides of the delta shock are incoming.

It follows from (\ref{2.19}) that
\begin{equation}\label{2.22}
\DF{dw(t)}{dt}=v_{\delta}(\rho_{+}-\rho_{-})-(\rho_{+}u_{+}-\rho_{-}u_{-}),
\end{equation}
\begin{equation}\label{2.23}
v_{\delta}\DF{dw(t)}{dt}=v_{\delta}(\rho_{+}u_{+}-\rho_{-}u_{-})-(\rho_{+}u_{+}^{2}-\rho_{-}u_{-}^{2})+A(u_{+}-u_{-}),
\end{equation}
Thus, we have
\begin{equation}\label{2.24}
(\rho_{+}-\rho_{-})v_{\delta}^{2}-2(\rho_{+}u_{+}-\rho_{-}u_{-})v_{\delta}+(\rho_{+}u_{+}^{2}-\rho_{-}u_{-}^{2})-A(u_{+}-u_{-})=0,
\end{equation}

For convenience, let us denote
\begin{equation}\label{2.25}
w_{0}=\sqrt{\rho_{+}\rho_{-}(u_{+}-u_{-})((u_{+}-u_{-})-(\frac{A}{\rho}_{+}-\frac{A}{\rho}_{-}))},
\end{equation}

If $\rho_+\neq\rho_-$, with the entropy condition (\ref{2.20}) in
mind, one can obtain directly from (\ref{2.24}) that
\begin{equation}\label{2.26}
v_{\delta}=\frac{\rho_+ u_+-\rho_-u_-+w_{0}}{\rho_+-\rho_-},
\end{equation}
which enables us to get
\begin{equation}\label{2.27}
\sigma(t)=v_{\delta}+\beta t,\ \ x(t)=v_{\delta}t+\frac{1}{2}\beta
t^{2}\ \ w(t)=w_{0}t,
\end{equation}

Otherwise,if $\rho_+=\rho_-$, then we have

\begin{equation}\label{2.28}
v_{\delta}=\frac{1}{2}(u_++u_--\frac{A}{\rho_{-}}).
\end{equation}
In this particular case, we can also get
\begin{equation}\label{2.29}
\sigma(t)=\frac{1}{2}(u_++u_--\frac{A}{\rho_{-}})+\beta t,\
x(t)=\frac{1}{2}(u_++u_--\frac{A}{\rho_{-}})t+\frac{1}{2}\beta
t^{2},\  w(t)=(\rho_-u_--\rho_+ u_+)t.
\end{equation}

\section{Riemann problem for the original system}
In this section, let us return to the Riemann problem (\ref{1.1})
and (\ref{1.2}). If $(\rho_+,u_+)\in$ I $\cup$ I\!I, namely
$u_{+}>u_{-}-\frac{A}{\rho_{-}}$, the Riemann solutions to
(\ref{1.1}) and (\ref{1.2}) can be represented as
\begin{equation}\label{3.1}
(\rho,u)(x,t)=\left\{\begin{array}{ll}
(\rho_-,u_-+\beta t),\ \ \ \ \ \ \ \ \ \ \ x<x_{1}(t),\\
(\rho_*,v_*+\beta t),\ \ \  \ \ \ \ \ \ \ \ x_{1}(t)<x<x_{2}(t),\\
(\rho_+,u_++\beta t),\ \ \ \ \ \ \ \ \ \ \ x>x_{2}(t),
\end{array}\right.
\end{equation}
where $(\rho_*,v_*)$ is given by (\ref{2.9}) and  the position of
the two contact discontinuities $J_{1}$ and $J_{2}$ are given
respectively by
\begin{equation}\label{3.2}
x_{1}(t)=(u_{-}-\frac{A}{\rho}_{-})t+\frac{1}{2}\beta t^{2},\ \
x_{2}(t)=u_{+}t+\frac{1}{2}\beta t^{2}.
\end{equation}
Let us draw Fig.3 to illustrate this situation in detail.

Analogously, if $(\rho_+,u_+)\in$ I\!I\!I, namely $0\leq u_{+}\leq
u_{-}-\frac{A}{\rho_{-}}$, then we can also define the weak
solutions in the sense of distributions to the Riemann problem
(\ref{1.1}) and (\ref{1.2}) below.

\unitlength 1mm 
\linethickness{0.4pt}
\ifx\plotpoint\undefined\newsavebox{\plotpoint}\fi 
\begin{picture}(169,69.25)(10,0)
\put(82,18){\vector(1,0){.07}}
\put(15,18){\line(1,0){67}}
\put(9,68){\vector(0,1){.07}}
\put(9,20.25){\line(0,1){47.75}}
\qbezier(40,18)(49.75,47.63)(75.5,57.25)
\qbezier(40,18)(16.13,34.75)(22.5,65.75)
\put(166.25,18.25){\vector(1,0){.07}}
\put(96,18.5){\line(1,0){68.78125}}
\put(91.25,68){\vector(0,1){.07}}
\put(91.25,21){\line(0,1){47}}
\qbezier(125,18.5)(126.63,54.25)(144.75,68)
\qbezier(125,18.5)(138.75,44.63)(161.25,55.5)
\put(5.5,65.25){$t$}
\put(84,18){$x$}
\put(85.75,68.75){$t$}
\put(169,18){$x$}
\put(23.5,69.25){$J_1$}
\put(77,59.5){$J_2$}
\put(39.75,14.75){0}
\put(125,15){0}
\put(136.25,66.75){$J_1$}
\put(164.5,56){$J_2$}
\put(10,23){$(\rho_-,u_-+\beta t)$}
\put(53.25,32.75){$(\rho_+,u_++\beta t)$}
\put(28,40){$(\rho_*,v_*+\beta t)$}
\put(100,40){$(\rho_-,u_-+\beta t)$}
\put(134,51.75){$(\rho_*,v_*+\beta t)$}
\put(145,32){$(\rho_+,u_++\beta t)$}
\put(45,12.75){\makebox(0,0)[cc]{(a)$\ \ u_--\frac{A}{\rho_-}<0<u_+$}}
\put(130,12.75){\makebox(0,0)[cc]{(b) $\ \ 0<u_--\frac{A}{\rho_-}<u_+$}}
\put(90,6){\makebox(0,0)[cc]
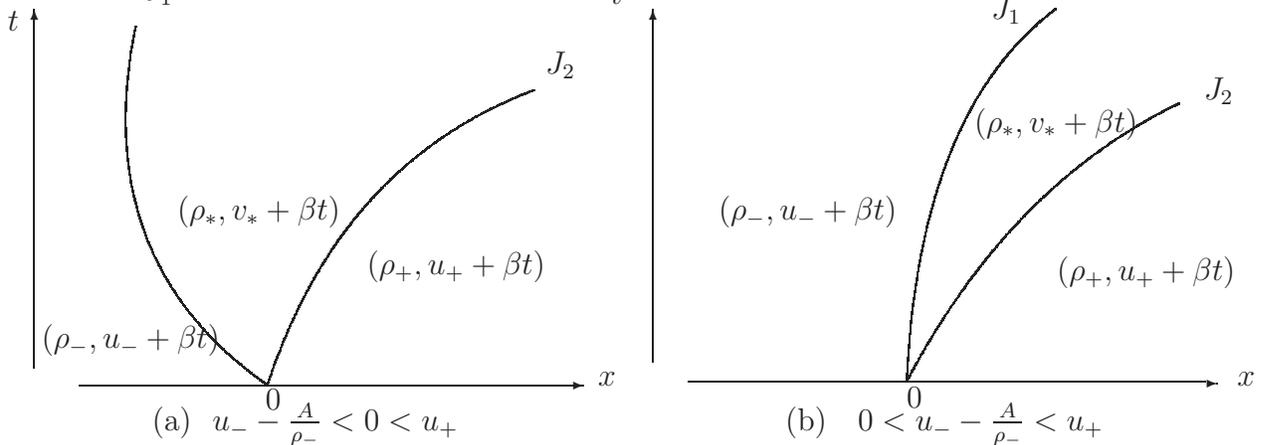
{Fig.3 The Riemann solution to (1.1) and (1.2) when $u_--\frac{A}{\rho_-}<u_+$ and $\beta>0$,}}
\put(59,3){\makebox(0,0)[cc]
{where $(\rho_*,v_*)$ is given by (2.9).}}

\end{picture}

\begin{defn}\label{defn:3.1}
Let $(\rho,u)$ be a pair of distributions in which $\rho$ has the
form of (\ref{2.11}), then it is called as the delta shock wave
solution to the Riemann problem (\ref{1.1}) and (\ref{1.2}) if it
satisfies
\begin{equation}\label{3.3}
\left\{\begin{array}{ll}
\langle\rho,\psi_{t}\rangle+\langle\rho u,\psi_{x}\rangle=0,\\
\langle\rho( u+P)),\psi_{t}\rangle+\langle\rho u(
u+P)),\psi_{x}\rangle=-\langle\beta\rho,\psi\rangle,
 \end{array}\right.
\end{equation}
for any $\psi\in C_0^\infty(R\times R^{+})$, in which
 $$\langle\rho
u(
u+P)),\psi\rangle=\int_{0}^{\infty}\int_{-\infty}^{\infty}(\widehat{\rho}u(u-\frac{1}{\widehat{\rho}}))\psi
dxdt+\langle w(t)(u_{\delta}(t))^{2}\delta_{S}\psi\rangle,$$ and
$u_{\delta}(t)$ is the assignment of $u$ on this delta shock wave
curve. \end{defn}

With the above definition in mind, if $u_{+}<
u_{-}-\frac{A}{\rho_{-}}$ is satisfied, then we look for a piecewise
smooth solution to the Riemann problem (\ref{1.1}) and (\ref{1.2})
in the form
\begin{equation}\label{3.4}
(\rho,u)(x,t)=\left\{\begin{array}{ll}
(\rho_-,u_-+\beta t),\ \ \ \ \ \ \ \ \ \ \ \ \ \ \ \ \ x<x(t),\\
(w(t)\delta(x-x(t)),u_{\delta}(t)),\ \ \  \ x=x(t),\\
(\rho_+,u_++\beta t),\ \ \ \ \ \ \ \ \ \ \ \ \ \ \ \ \ x>x(t),
\end{array}\right.
\end{equation}
It is worthwhile to notice that $u_{\delta}(t)-\beta t$ is assumed
to be a constant based on the result in Sect.2. With the similar
analysis and derivation as before, the delta shock wave solution of
the form (\ref{3.4}) to the Riemann problem (\ref{1.1}) and
(\ref{1.2}) should also satisfy the following generalized
Rankine-Hugoniot conditions
\begin{equation}\label{3.5}
\left\{\begin{array}{ll}
\DF{dx(t)}{dt}=\sigma(t)=u_{\delta}(t),\\[4pt]
\DF{dw(t)}{dt}=\sigma(t)[\rho]-[\rho u],\\[4pt]
\DF{d(w(t)u_{\delta}(t))}{dt}=\sigma(t)[\rho(
u-\frac{A}{\rho})]-[\rho u( u-\frac{A}{\rho}) ]+\beta w(t).
\end{array}\right.
\end{equation}
in which the jumps across the discontinuity are
\begin{equation}\label{3.6}
[\rho u]=\rho_+(u_++\beta
t)-\rho_-(u_-+\beta t),
\end{equation}
\begin{equation}\label{3.7}
[\rho u( u-\frac{A}{\rho})]=\rho_+(u_++\beta t)(u_++\beta
t-\frac{A}{\rho_{+}})-\rho_-(u_-+\beta t)(u_-+\beta
t-\frac{A}{\rho}).
\end{equation}

In order to ensure the uniqueness to the Riemann problem (\ref{1.1})
and (\ref{1.2}), the over-compressive entropy condition for the
delta shock wave
\begin{equation}\label{3.8}
u_{+}+\beta t<u_{\delta}(t)<u_{-}-\frac{A}{\rho_{-}}+\beta t.
\end{equation}
should also be proposed when $0\leq u_{+}<
u_{-}-\frac{A}{\rho_{-}}$.

Like as before, we can also obtain $x(t),\sigma(t)$ and $w(t)$ from
(\ref{3.5}) and (\ref{3.8}) together. In brief, we have the
following theorem to depict the Riemann solution to (\ref{1.1}) and
(\ref{1.2}) when the Riemann initial data (\ref{1.2}) satisfy $0\leq
u_{+}< u_{-}-\frac{A}{\rho_{-}}$ and $\rho_{+}\neq\rho_{-}$.

\begin{thm}\label{thm:3.1} If both $0\leq u_{+}< u_{-}-\frac{A}{\rho_{-}}$ and
$\rho_{+}\neq\rho_{-}$ are satisfied, then the delta shock solution
to the Riemann solutions to (\ref{1.1}) and (\ref{1.2}) can be
expressed as
\begin{equation}\label{3.9}
\left\{\begin{array}{ll}
\DF{dx(t)}{dt}=\sigma(t)=u_{\delta}(t),\\[4pt]
\DF{dw(t)}{dt}=\sigma(t)[\rho]-[\rho u],\\[4pt]
\DF{d(w(t)u_{\delta}(t))}{dt}=\sigma(t)[\rho(
u-\frac{A}{\rho})]-[\rho u( u-\frac{A}{\rho}) ]+\beta w(t).
\end{array}\right.
\end{equation}
in which
\begin{equation}\label{3.10}
\sigma(t)=u_{\delta}(t)=v_{\delta}+\beta t,\ \
x(t)=v_{\delta}t+\frac{1}{2}\beta t^{2}\ \ w(t)=w_{0}t,
\end{equation}
in which $w_{0}$ and $v_{\delta}$ are given by (\ref{2.21}) and
(\ref{2.22}) respectively. \end{thm}

Let us check briefly that the above constructed delta shock wave
solution (\ref{3.9}) and (\ref{3.10}) should satisfy (\ref{1.1}) in
the sense of distributions. The proof of this theorem is completely
analogs to those in \cite {Shen1,Shen2}. Therefore, we only deliver
the main steps for the proof of the second equality in (\ref{3.3})
for completeness. Actually, one can deduce that

\begin{eqnarray}
I&=&\int_{0}^{\infty}\int_{-\infty}^{\infty}(\rho(
u-\frac{A}{\rho})\psi_{t}+\rho u(
u-\frac{A}{\rho}) \psi_{x})dxdt\nonumber\\
&=&\int_{0}^{\infty}\int_{-\infty}^{x(t)}(\rho_{-}( u_{-}+\beta
t-\frac{A}{\rho_{-}})\psi_{t}+\rho_{-}(u_{-}+\beta t)(u _{-}+\beta
t-\frac{A}{\rho_{-}})\psi_{x})dx dt \nonumber\\
&&+\int_{0}^{\infty}\int^{\infty}_{x(t)}(\rho_{+}( u_{+}+\beta
t-\frac{A}{\rho_{+}})\psi_{t}+\rho_{+}(u_{+}+\beta t)(u _{+}+\beta
t-\frac{A}{\rho_{+}}) \psi_{x})dx dt \nonumber\\
&&+\int_{0}^{\infty}w_{0}t(v_{\delta}+\beta
t)(\psi_{t}(x(t),t)+(v_{\delta}+\beta
t)\psi_{x}(x(t),t))dt\nonumber.
\end{eqnarray}

It can be derived from (\ref{3.10}) that the curve of delta shock
wave is given by
\begin{equation}\label{3.11}
x(t)=v_{\delta}t+\frac{1}{2}\beta t^{2}.
\end{equation}

\unitlength 1mm 
\linethickness{0.5pt}
\ifx\plotpoint\undefined\newsavebox{\plotpoint}\fi 
\begin{picture}(182,63)(10,0)
\put(85.5,12.5){\vector(1,0){.07}}
\put(12,12.5){\line(1,0){73.5}}
\put(4.75,61){\vector(0,1){.07}}
\put(4.75,15){\line(0,1){46}}
\qbezier(38.25,12.5)(39.75,40.25)(67.25,55.25)
\put(176,13.5){\vector(1,0){.07}}
\put(104,13.5){\line(1,0){72}}
\put(94.75,63){\vector(0,1){.07}}
\put(94.75,14.5){\line(0,1){48.5}}
\put(1.75,58.75){$t$}
\put(87.75,11.75){$x$}
\put(66,55){$\delta \!S$}
\put(92,61.75){$t$}
\put(177,13){$x$}
\put(110,50){$\delta \!S$}
\put(21,35.75){$(\rho_-,u_-+\beta t)$}
\put(55,35.75){$(\rho_+,u_++\beta t)$}
\put(106.75,35){$(\rho_-,u_-+\beta t)$}
\put(150,35.25){$(\rho_+,u_++\beta t)$}
\put(38,8.5){0}
\put(115,9){0}
\put(45,8.5){(a)$\ \ \beta>0$}
\put(127.5,8.5){(b)$\ \ \beta<0$}
\put(90,5){\makebox(0,0)[cc]
{Fig.4 The delta shock wave solution to (1.1) and (1.2) when $u_+<u_--\frac{A}{\rho_-}$, }}
 \put(90,1){\makebox(0,0)[cc]
{ where $v_\delta>0$ is given by (2.26) for $\rho_-\neq\rho_+$ and (2.28) for $\rho_-= \rho_+$.}}
 \qbezier(115,13.5)(182,22.75)(105.5,57.75)
\end{picture}

For $\beta>0$ (Fig.4a), there exists an inverse
function of $x(t)$ globally in the time $t$, which may be written in
the form
$$t(x)=\sqrt{\frac{v_{\delta}^{2}}{\beta^{2}}+\frac{2x}{\beta}}-\frac{v_{\delta}}{\beta}.$$
Otherwise, for $\beta<0$ (Fig.4b), there is a critical
point $(-\frac{v_{\delta}^{2}}{2\beta},-\frac{v_{\delta}}{\beta})$
on the delta shock wave curve such that $x'(t)$ change its sign when
across the critical point. Thus, the inverse function of $x(t)$ is
needed to find respectively for $t\leq-\frac{v_{\delta}}{\beta}$ and
$t>-\frac{v_{\delta}}{\beta}$, which enable us to have
\begin{equation}\nonumber
t(x)=\left\{\begin{array}{ll}
-\sqrt{\frac{v_{\delta}^{2}}{\beta^{2}}+\frac{2x}{\beta}}-\frac{v_{\delta}}{\beta},\ \ \ t\leq-\frac{v_{\delta}}{\beta},\\[4pt]
\sqrt{\frac{v_{\delta}^{2}}{\beta^{2}}+\frac{2x}{\beta}}-\frac{v_{\delta}}{\beta},\
\  \ \ \ t>-\frac{v_{\delta}}{\beta}.
\end{array}\right.
\end{equation}

Without loss of generality, let us assume that $\beta>0$ for
simplicity. Actually, the other situation can be dealt with
similarly. Under our assumption, it follows from (\ref{3.12}) that
the position of delta shock wave satisfies $x=x(t)>0$ for all the
time. It follows from (\ref{3.10}) that
\begin{eqnarray}
\frac{d\psi(x(t),t)}{dt}&=&\psi_t(x(t),t)+\frac{dx(t)}{dt}\psi_x(x(t),t)\nonumber\\
&=&\psi_t(x(t),t)+(v_\delta+\beta t)\psi_x(x(t),t)\nonumber\\
&=&\psi_t(x(t),t)+u_\delta(t)\psi_x(x(t),t)\nonumber.
\end{eqnarray}

By exchanging the ordering of integrals and using integration by
parts, we have

\begin{eqnarray}
I&=&\int_{0}^{\infty}\int^{\infty}_{t(x)}\rho_{-}( u_{-}+\beta
t-\frac{A}{\rho_{-}})\psi_{t}dt dx
+\int_{0}^{\infty}\int^{\infty}_{t(x)}\rho_{-}(u_{-}+\beta t)(u
_{-}+\beta t-\frac{A}{\rho_{-}})\psi_{x}dt dx \nonumber\\
&&+\int_{0}^{\infty}\int_{0}^{t(x)}\rho_{+}( u_{+}+\beta
t-\frac{A}{\rho_{+}})\psi_{t}dt
dx+\int_{0}^{\infty}\int_{0}^{t(x)}\rho_{+}(u_{+}+\beta t)(u
_{+}+\beta t-\frac{A}{\rho_{+}}) \psi_{x}dtdx\nonumber\\
&&+\int_{0}^{\infty}w_{0}t(v_{\delta}+\beta
t)d\psi(x(t),t)\nonumber\\
&=&\int_{0}^{\infty}(\rho_{+}( u_{+}+\beta
t(x)-\frac{A}{\rho_{+}})-\rho_{-}( u_{-}+\beta
t(x)-\frac{A}{\rho_{-}}))\psi(x,t(x))dx\nonumber\\
&&+\int_{0}^{\infty}(\rho_{-}(u_{-}+\beta t)( u_{-}+\beta
t-\frac{A}{\rho_{-}})-\rho_{+}(u_{+}+\beta t)( u_{+}+\beta
t-\frac{A}{\rho_{+}}))\psi(x(t),t)dt \nonumber\\
&&-\int_{0}^{\infty}\int^{\infty}_{t(x)}\beta\rho_{-}\psi dt
dx-\int_{0}^{\infty}\int_{0}^{t(x)}\beta\rho_{+}\psi dt
dx-\int_{0}^{\infty}w_{0}(v_{\delta}+2\beta
t)\psi(x(t),t)dt\nonumber\\
&=&\int_{0}^{\infty}A(t)\psi(x(t),t)dt-\beta(\int_{0}^{\infty}\int_{-\infty}^{x(t)}\rho_{-}\psi
dx dt+\int_{0}^{\infty}\int^{\infty}_{x(t)}\rho_{+}\psi dx
dt)\label{3.12},
\end{eqnarray}
in which \begin{eqnarray}C(t)&=&(\rho_{+}( u_{+}+\beta
t-\frac{A}{\rho_{+}})-\rho_{-}( u_{-}+\beta
t-\frac{A}{\rho_{-}}))(v_{\delta}+\beta t)\nonumber\\
&&+(\rho_{-}(u_{-}+\beta t)( u_{-}+\beta
t-\frac{A}{\rho_{-}})-\rho_{+}(u_{+}+\beta t)( u_{+}+\beta
t-\frac{A}{\rho_{+}}))\nonumber\\
&&-w_{0}(v_{\delta}+2\beta t)\nonumber.
\end{eqnarray}
By a tedious calculation, we have
\begin{equation}\label{3.13}
A(t)=-\beta w_{0}t=-\beta w(t).
\end{equation}
Thus, it can be concluded from (\ref{3.12}) and (\ref{3.13})
together that the second equality in (\ref{3.3}) holds in the sense
of distributions. The proof is completed.

\begin{rem}\label{rem:3.1}
If both $0\leq u_{+}<
u_{-}-\frac{A}{\rho_{-}}$ and $\rho_{+}=\rho_{-}$ are satisfied,
then the delta shock solution to the Riemann solutions to
(\ref{1.1}) and (\ref{1.2}) can be expressed in the form
(\ref{3.11}) where
\begin{equation}\label{3.14}
\sigma(t)=u_{\delta}(t)=\frac{1}{2}(u_++u_--\frac{A}{\rho_{-}})+\beta
t,\ x(t)=\frac{1}{2}(u_++u_--\frac{A}{\rho_{-}})t+\frac{1}{2}\beta
t^{2},\ w(t)=(\rho_-u_--\rho_+ u_+)t.
\end{equation}
The process of proof is completely similar and we omit the details.
\end{rem}

\begin{rem}\label{rem:3.2}   If $u_{+}=u_{-}-\frac{A}{\rho_{-}}$, then the
delta shock solution to the Riemann solutions to (\ref{1.1}) and
(\ref{1.2}) can be also expressed as the form in Theorem \ref{thm:3.1} and Remark \ref{rem:3.1}.                                                                                                                                           The process of proof is easy and we
omit the details.
\end{rem}

\section{The vanishing pressure limit of Riemann solutions to (\ref{1.1}) and (\ref{1.2})}

In this section, we consider the vanishing pressure limit of Riemann
solutions to (\ref{1.1}) and (\ref{1.2}). According to the relations
between $u_-$ and $u_+$, we will divide our discussion into the
following three case:

 (1) $u_-<u_+$; \ \ \ \ \ \     (2) $u_-=u_+$;  \ \ \ \  \ \    (3) $u_->u_+$.

\textbf{Case 4.1.}  $u_-<u_+$

In this case, $(\rho_+,u_+)\in$ I in the $(\rho,v)$ plane, so the
Riemann solutions to (\ref{1.1}) and (\ref{1.2}) is given by
(\ref{3.1}) and (\ref{3.2}), where $(\rho_*,v_*)$ is given by
(\ref{2.9}). From (\ref{2.9}) we have

$$\lim\limits_{A\rightarrow0}\rho_*=\lim\limits_{A\rightarrow0}\frac{A}{u_+-u_-+\frac{A}{\rho_-}}=0,$$
which indicates the occurrence of the vacuum states. Furthermore,
the Riemann solutions to (\ref{1.1}) and (\ref{1.2}) converge to

\begin{equation}\label{4.1}
\lim\limits_{A\rightarrow0}(\rho,u)(x,t)=\left\{\begin{array}{ll}
(\rho_-,u_-+\beta t),\ \ \ \ \ \ \ \ \ \ \ x<u_{-}t+\frac{1}{2}\beta t^{2},\\
vacuum,\ \ \ \ \ \ \ \ \ \ \ \ \ \ \ \ \ u_{-}t+\frac{1}{2}\beta t^{2}<x<u_{+}t+\frac{1}{2}\beta t^{2},\\
(\rho_+,u_++\beta t),\ \ \ \ \ \ \ \ \ \ \ x>u_{+}t+\frac{1}{2}\beta
t^{2},
\end{array}\right.
\end{equation}
which is exactly the corresponding Riemann solutions to the
transportation equations with the same source term and the same
initial data \cite{Shen1}.

\textbf{Case 4.2.}  $u_-=u_+$

In this case, $(\rho_+,u_+)$ is on the $J_{2}$ curve in the
$(\rho,v)$ plane, so the Riemann solutions to (\ref{1.1}) and
(\ref{1.2}) is given as follows:
\begin{equation}\label{4.2}
(\rho,u)(x,t)=\left\{\begin{array}{ll}
(\rho_-,u_-+\beta t),\ \ \ \ \ \ \ \ \ \ \ x<u_{-}t+\frac{1}{2}\beta t^{2},\\
(\rho_+,u_++\beta t),\ \ \ \ \ \ \ \ \ \ \ x>u_{+}t+\frac{1}{2}\beta
t^{2},
\end{array}\right.
\end{equation}
which is exactly the corresponding Riemann solutions to the
transportation equations with the same source term and the same
initial data \cite{Shen1}.

\textbf{Case 4.3.}  $u_->u_+$
\begin{lem}\label{lem:4.1}
If $u_->u_ +$, there exist $A_{1}>A_{0}>0$, such that
$(\rho_+,u_+)\in$ I\!I as $A_{0}<A<A_{1}$, and $(\rho_+,u_+)\in$
I\!I\!I as $A\leq A_{0}$.
\end{lem}

\noindent\textbf{Proof.} If $(\rho_+,u_+)\in$ I\!I , then
$0<u_{-}-\frac{A}{\rho_{-}}<u_+<u_-$, which gives
$\rho_{-}u_{-}>A>\rho_{-}(u_--u_+)$. Thus we take
$A_0=\rho_{-}(u_--u_+)$ and $A_0=\rho_{-}u_-$, then
$(\rho_+,u_+)\in$ I\!I as $A_{0}<A<A_{1}$ and $(\rho_+,u_+)\in$
I\!I\!I as $A\leq A_{0}$.

When $A_{0}<A<A_{1}$, $(\rho_+,u_+)\in$ I\!I in the $(\rho,v)$
plane, so the Riemann solutions to (\ref{1.1}) and (\ref{1.2}) is
given by (\ref{3.1}) and (\ref{3.2}), where $(\rho_*,v_*)$ is given
by (\ref{2.9}). From (\ref{2.9}) we have From (\ref{2.9}) we have

$$\lim\limits_{A\rightarrow A_0}\rho_*=\lim\limits_{A\rightarrow A_0}\frac{A}{u_+-u_-+\frac{A}{\rho_-}}
=\lim\limits_{A\rightarrow A_0}\frac{\rho_-A}{A-A_{0}}=\infty.$$

Furthermore, we have the following result.
\begin{lem}\label{lem:4.2}
Let $\frac{dx_1(t)}{dt}=\sigma_1(t)$,
$\frac{dx_2(t)}{dt}=\sigma_2(t)$, then we have
\begin{equation}\label{4.3}
\lim\limits_{A\rightarrow A_0}v_*+\beta t=\lim\limits_{A\rightarrow
A_0}\sigma_1(t)=\lim\limits_{A\rightarrow
A_0}\sigma_2(t)=(u_{-}-\frac{A_0}{\rho_{-}})t+\beta t=u_++\beta
t=:\sigma(t),
\end{equation}
\begin{equation}\label{4.4}
\lim\limits_{A\rightarrow
A_0}\int_{x_1(t)}^{x_2(t)}\rho_*dx=A_0 t,
\end{equation}
\begin{equation}\label{4.5}
\lim\limits_{A\rightarrow A_0}\int_{x_1(t)}^{x_2(t)}\rho_*(v_*+\beta
t)dx=(u_++\beta t)A_0 t.
\end{equation}
\end{lem}

\noindent\textbf{Proof.} (\ref{4.3}) is obviously true. We will only
prove (\ref{4.4}) and (\ref{4.5}).

\begin{equation}\nonumber
\lim\limits_{A\rightarrow
A_0}\int_{x_1(t)}^{x_2(t)}\rho_*dx=\lim\limits_{A\rightarrow
A_0}\rho_*(x_2(t)-x_1(t))=\lim\limits_{A\rightarrow
A_0}\frac{A}{u_+-u_-+\frac{A}{\rho_-}}(u_+-u_-+\frac{A}{\rho_-})t=A_0
t,
\end{equation}
\begin{eqnarray}
\lim\limits_{A\rightarrow A_0}\int_{x_1(t)}^{x_2(t)}\rho_*(v_*+\beta
t)dx=(u_++\beta t)\lim\limits_{A\rightarrow
A_0}\int_{x_1(t)}^{x_2(t)}\rho_*dx=(u_++\beta t)A_0 t\nonumber.
\end{eqnarray}
The proof is completed.

It can be concluded from Lemma \ref{lem:4.2} that the curves of the
two contact discontinuities $J_1$ and $J_2$ will coincide when $A$
tends to $A_0$ and the delta shock waves will form.  Next we will
arrange the values which gives the exact position, propagation speed
and strength of the delta shock wave according to  Lemma
\ref{lem:4.2}.

 From (\ref{4.4}) and (\ref{4.5}), we let
\begin{equation}\label{4.6}
w(t)=A_0 t,
\end{equation}
\begin{equation}\label{4.7}
w(t)u_\delta(t)=(u_++\beta t)A_0 t,
\end{equation}
then
\begin{equation}\label{4.8}
u_\delta(t)=(u_++\beta t),
\end{equation}
which is equal to $\sigma(t)$. Furthermore, by letting
$\frac{dx(t)}{dt}=\sigma(t)$, we have
\begin{equation}\label{4.9}
x(t)=u_+ t+\frac{1}{2}\beta t^2.
\end{equation}

From (\ref{4.6})-(\ref{4.9}), we can see that the quantities defined
above are exactly consistent with those given by
(\ref{2.25})-(\ref{2.29}) or (\ref{3.10}) in which we take  $A=A_0$.
Thus, it uniquely determines that the limits of the Riemann
solutions to the system (\ref{1.1}) and (\ref{1.2}) when
$A\rightarrow A_0$ in the case $(\rho_+,u_+)\in$ I\!I is just the
delta shock solution of (\ref{1.1}) and (\ref{1.2}) in the case
$(\rho_+,u_+)\in$ S, where S is actually the boundary between the
regions I\!I and I\!I\!I. So we get the following results in the
case $u_+<u_-$.

\begin{thm}\label{thm:4.1} If $u_{+}< u_{-}$, for each fixed $A$ with $A_{0}<A<A_{1}$, $(\rho_+,u_+)\in$
\rm{I\!I}
assuming that $(\rho,u)$ is a solution containing two contact
discontinuities $J_1$ and $J_2$ of (\ref{1.1}) and (\ref{1.2}) which
is constructed in Section 3, it is obtained that when $A\rightarrow
A_0$, $(\rho,u)$ converges to a delta shock wave solution of
(\ref{1.1}) and (\ref{1.2}) when $A=A_0$.
\end{thm}

When $A\leq A_{0}$, $(\rho_+,u_+)\in$ I\!I\!I, so the Riemann
solutions to (\ref{1.1}) and (\ref{1.2}) is given by (\ref{3.4})
with (\ref{3.10}) or (\ref{3.14}), which is a delta shock wave
solution. It is easy to see that when $A\rightarrow 0$, for
$\rho_+\neq \rho_-$,
$$x(t)\rightarrow \sigma t+\frac{1}{2}\beta
t^{2},\ \ w(t)\rightarrow\sqrt{\rho_+\rho_-}(u_-- u_+)t,\ \
\sigma(t)=u_{\delta}(t)\rightarrow  \sigma +\beta t,$$ where
 $\sigma=\frac{\sqrt{\rho_-}u_-+\sqrt{\rho_+}u_+}{\sqrt{\rho_-}+\sqrt{\rho_+}}$,
for $\rho_+= \rho_-$,
$$x(t)\rightarrow \frac{1}{2}(u_++u_-)t+\frac{1}{2}\beta
t^{2},\ \ w(t)\rightarrow \rho_+(u_-- u_+)t,\ \
\sigma(t)=u_{\delta}(t)\rightarrow  \frac{1}{2}(u_++u_-)+\beta t,$$
which is exactly the corresponding Riemann solutions to the
transportation equations with the same source term and the same
initial data \cite{Shen1}. Thus, we have the following result:

\begin{thm}If $u_{+}< u_{-}$, for each fixed $A<A_0$, $(\rho_+,u_+)\in$
\rm{I\!I\!I}
assuming that $(\rho,u)$ is a  a delta shock wave solution of
(\ref{1.1}) and (\ref{1.2}) which is constructed in Section 3, it is
obtained that when $A\rightarrow 0$, $(\rho,u)$ converges to a delta
shock wave solution to the transportation equations with the same
source term and the same initial data \cite{Shen1}.
\end{thm}

Now we summarize the main result in this section as follows.
\begin{thm}\label{thm:4.3}
As the pressure vanishes, the Riemann solutions to the
Chaplygin pressure Aw-Rascle model with Coulomb-like friction tend
to the three kinds of Riemann solutions to the transportation
equations with the same source term and the same initial data, which
included a delta shock wave and a vacuum state.
\end{thm}

\section{Conclusions and Discussions}
In this work, we have considered the solutions of the Riemann
problem for the Chaplygin pressure Aw-Rascle model with Coulomb-like
friction in the fully explicit form. In particular, the delta shock
wave solution has been discovered in some certain situations, which
may be used to explain the serious traffic jam. We find that the
Coulomb-like friction term takes the effect to curve the
characteristics to the parabolic curves such that the delta shock
wave discontinuity is also curved. Thus, the Riemann solutions
(\ref{1.1}) and (\ref{1.2}) are not self-similar any more. It is
worthwhile to note that the Riemann solutions of (\ref{1.1}) and
(\ref{1.2}) converge to the corresponding ones of the Chaplygin
pressure Aw-Rascle model as $\beta\rightarrow0$, namely the
Coulomb-like friction term vanishes. Finally, we analyze the
formation of $\delta$-shocks and vacuum states in the Riemann
solutions in the vanishing pressure limit and show that the Riemann
solutions of (\ref{1.1}) and (\ref{1.2}) converge to the
corresponding ones of the transportation equations with the same
source term as the pressure vanishes. These results generalize those
obtained in \cite{Chen-Liu1,Sheng-Wang-Yin} for homogeneous
equations to nonhomogeneous equations and are also applicable to the nonsymmetric system of
Keyfitz-Kranzer type with the same Chaplygin pressure and Coulomb-like friction.

It is interesting to notice that the above constructed Riemann
solutions of (\ref{1.1}) and (\ref{1.2}) can be obtained directly
from the ones of the Riemann problem for the homogeneous situation
by using the change of variables $x\rightarrow x-\frac{1}{2}\beta
t^{2}$ and $u\rightarrow u-\beta t$ together, see
\cite{Karelsky-Petrosyan-Tarasevich}. It also pointed out in
\cite{Karelsky-Petrosyan-Tarasevich} that these solutions are
drastically different from each other in that the characteristics
are the parabolic curves for the inhomogeneous situation instead of
the straight lines for the homogeneous situation. Furthermore, the
regions of constant flow are transformed into the regions of
constantly accelerated flow and the contact discontinuities and the
the delta shock waves bend into parabolic shapes under the influence
of the Coulomb-like friction term.

It is worthwhile to note that the method developed in this paper can
be used to the inhomogeneous Aw-Rascle model with generalized
Chaplygin pressure. Especially, the Aw-Rascle model with generalized
Chaplygin pressure has a significant mathematical difference with
the Aw-Rascle model with Chaplygin pressure. Thus, it is interesting
to study the Riemann problem for the Aw-Rascle model with
generalized Chaplygin pressure under the influence of the the
Coulomb-like friction term, whose resuts will
also be applicable to the nonsymmetric system of Keyfitz-Kranzer type (\ref{1.4}) with the same
pressure and Coulomb-like friction. We leave this problem for our future
work.

%
%
%
%







\end{document}